\documentclass[12pt]{article}
\usepackage{amsmath,amsfonts,amssymb,amsthm,amscd}
\title{Remarks on Tao's algebraic regularity lemma}
\date{October 28 2013}
\author{Anand Pillay\\University of Notre Dame \and Sergei Starchenko\\University of Notre Dame}
\newtheorem{Theorem}{Theorem}[section]

\newtheorem{Remark}[Theorem]{Remark}
\newtheorem{Lemma}[Theorem]{Lemma}
\newtheorem{Corollary}[Theorem]{Corollary}

\begin{document}
\maketitle

\begin{abstract} 
We give a stability-theoretic proof of the algebraic regularity lemma from \cite{Tao}, in a slightly strengthened form.  We also point out that the underlying lemmas hold at a greater level of generality, namely ``measurable" theories and structures in the sense of Elwes-Macpherson-Steinhorn. 

\end{abstract}

\section{Introduction}

If $K$ is a pseudofinite field (infinite model of the theory of finite fields) then we can attach to any definable set $X$ a dimension $d$ and ``measure" $\mu$, which come from counting points in  finite fields \cite{CDM}. Pseudofinite fields have simple first order theory, and they are also geometric structures in the sense of \cite{HP} and the corresponding notions of independence coincide.  This was all made very clear in the 1990's but we recall some of these facts later in the introduction. We now present the results.

\begin{Lemma}  Let $F$ be a saturated pseudofinite field, and $A$ a small algebraically closed set. Let $\phi_{0}(x)$ be a formula over $A$ with dimension $n$, and let $\phi(x,y)$, $\psi(x,z)$ be formulas over $A$ each of which imply $\phi_{0}$.  Let $p(y)$ and $q(z)$ be complete types over $A$. Suppose that for some independent realizations $a$ of $p$ and $b$ of $q$, the (dimension, measure) of $\phi(x,a)\wedge\psi(x,b)$ is $(n,r)$. Then for all independent realizations $a$ of $p$ and $b$ of $q$, the (dimension, measure) of $\phi(x,a)\wedge\psi(x,b)$ equal to $(n,r)$. 
\end{Lemma}

The next corollary is  Proposition 27 of \cite{Tao} (first reduction of the regularity lemma), but for pseudofinite fields of any characteristic and also with a control over parameters of definition.

\begin{Corollary} Let $K$ be a pseudofinite field.   Let $V,W$, and $E\subseteq V\times W$ be $A$-definable sets. Assume $dim(V) = n$ and $dim(W) = k$. Then we can partition $W$ into $acl(A)$-definable sets $W_{1},...,W_{m}$ such that for each $1\leq i,j\leq m$, there is 
$c_{i.j}>0$  and an $acl(A)$-definable subset $D_{i.j}$ of $W_{i}\times W_{j}$ with $dim(D_{i,j}) <2k$ such that either $dim(E(x,a)\cap E(x,b)) < n$ for all $(a,b)\in (W_{i}\times W_{j})\setminus D_{i,j}$ or
the (dimension, measure) of $(E(x,a)\cap E(x,b))$ equals $(n,c_{i,j})$ for all $(a,b)\in W_{i}\times W_{j}\setminus D_{i,j}$. 

\end{Corollary}

As in \cite{Tao}  we conclude the following (which improves Lemma 5 of \cite{Tao} by replacing ``$\bf F$ of characteristic at least $C$"  by ``$\bf F$ of cardinality at least $C$"). We could also give conditions on parameters of definition. See \cite{Tao} for the notation.

\begin{Corollary} If $M> 0$, there exists $C = C_{M} > 0$ such that: whenever $\bf F$ is a finite field of {\em cardinality} $\geq C$, $V, W$ are nonempty definable sets in $\bf F$ of complexity at most $M$ and $E\subseteq V\times W$ is another definable set of complexity at most $M$, then there exists partitions $V = V_{1}\cup ... \cup V_{a}$, $W = W_{1} \cup .. \cup W_{b}$ such that:
\newline
(i) For all $i=1,..,a$ and $j=1,..,b$, we have $|V_{i}| \geq |V|/|C|$ and $|W_{j}|\geq |W|/|C|$.
\newline
(ii) The $V_{i}$'s and $W_{j}$'s are definable with complexity at most $C$.
\newline
(iii) For all $i,j$, $A\subset V_{i}$ and $B\subset W_{j}$, $| |E\cap(A\times B)| - d_{ij}|A||B|| \leq C|{\bf F}|^{-1/4}|V_{i}||W_{j}|$,
where $d_{ij} = |E\cap(V_{i}\times W_{j})|/|V_{i}||W_{j}|$.

\end{Corollary}

\vspace{5mm}
\noindent
Recall the main theorem of \cite{CDM}. Here $L$ is the language of rings. Let $\phi(x,y)$ be an $L$-formula, where $x$ is an $n$-tuple of variables, and $y$ an arbitrary tuple. Then there are a positive constant $C$,  a finite set $D$ of pairs $(d,\mu)$ where $d$ is a nonnegative integer $\leq n$, and $\mu$ a positive rational, and for each $(d,\mu)\in D$ an $L$-formula $\psi_{d,\mu}(y)$, such that 
\newline
(i) the $\psi$'s partition $y$-space, 
\newline
(ii) for any finite field $\bf F$, $(d,\mu)\in D$, and tuple $b$ from $\bf F$ such that ${\bf F}\models \psi_{d,\mu}(b)$, we have that
$||\phi(x,b)({\bf F})| - \mu({\bf F}^{d})| \leq C {\bf F}^{d-(1/2)}$.

\vspace{2mm}
\noindent
So all this applies to a pseudofinite field, to give the invariants $(dim,measure)$ of any definable set. In fact the dimension of a definable set is the algebraic geometric dimension of its Zariski closure.  This dimension gives a notion of independence: $a$ independent from $b$ over $A$ if $dim(tp(a/A,b)) = dim(tp(a/A))$, which coincides with nonforking. (See section 5 of \cite{HP}.)  The dimension, measure function has a number of properties, including a Fubini statement for definable surjections $f:X\to Y$. Elwes, Macpherson and Steinhorm abstracted these properties to give the notion of a {\em measurable} structure or theory \cite{Elwes}. It will again be simple of finite $SU$-rank and one can take dimension to be $SU$-rank.

\section{Proofs}

The only things required are to prove Lemma 1.1 and deduce Corollary 1.2.  The proof of 1.1 uses ``local stability". The material needed for the case at hand is all in \cite{HP}, but also in \cite{Pillay-book}.  The key new fact is Proposition 2.25 from \cite{Udi-approximate}. 

\vspace{5mm}
\noindent
{\em Proof of Lemma 1.1.} 
\newline
For any formula $\chi(x)$ with parameters which implies $\phi_{0}(x)$, let $\mu^{*}(\chi(x)) = 0$ if $dim(\chi(x)) < n$ and let $\mu^{*}(\chi(x)) = r$ if the dimension, measure of $\chi(x)$ is $(n,r)$. 
Then $\mu^{*}$ is a Keisler measure on $\phi_{0}(x)$ and is  definable over $\emptyset$ (so also over $A$) so in particular is invariant.Let $\phi(x,y)$ and $\psi(x,z)$ be as in the statement of the Lemma. For a given  $r$, $\{(a,b): \mu^{*}(\phi(x,a)\wedge \psi(x,b)) = r\}$ is an $A$-definable set,  defined by formula $\delta(y,z)$ say. Proposition 2.25 of \cite{Udi-approximate} says that $\delta(y,z)$ is stable. There should be a direct proof of this in the case at hand but we did not try to find it yet. We now argue as in the proof of the Independence Theorem in Lemma 5.22 of \cite{HP}, making use of results in section 5 of that paper.  Suppose for some independent realizations $a$ and $b$ of $p$ and $q$ respectively, $\models\delta(a,b)$. Independence coincides with nonforking. As $A$ is algebraically closed,  $tp_{\delta}(a/A)$ has a unique nonforking extension to a complete $\delta$-type over $A,b$.  Hence for any realization $a'$ of $p$ such that $a'$ is independent from $b$ over $A$ we have that $\models\delta(a',b)$, which suffices to prove Lemma 1.1.

\begin{Remark} The same proof of the identical statement works for ``measurable" theories.  However we go outside the ``geometric structures" context so cannot appeal directly to \cite{HP}. 
We simply have to know that in measurable structure we can take the dimension function to be $SU$-rank, so dimension independence corresponds to nonforking etc.

\end{Remark}

\vspace{5mm}
\noindent
{\em Proof of Corollary 1.2.}
\newline
This is just a routine application of compactness using Lemma 1.1 and definability of dimension, measure.  But we give a few details.   First let $S$ be the set of complete types over $acl(A)$ extending ``$x\in W$".  For each $p\in S$, we have by Lemma 1.1, that either for all independent (over $A$) realizations $a,b$ of $p$, the dimension of $E(x,a)\cap E(x,b)$ is $< n$, OR for all independent (over $A$) realizations $a, b$ of $p$ that the dimension, measure  of $E(x,a)\cap E(x,b)$  equals $(n,r)$ for fixed positive rational $r$. 
Let's apply compactness to the second possibility, as an example:
We have  the implication $p(y)$ and $p(z)$ and $\{\neg(\chi(y,z)):\chi(y,z)$ over $acl(A)$ and $dim(\chi(y,z)) < 2k\}$ implies   $\delta(y,z)$, where $\delta(y,z)$ says that dimension, measure of $E(x,y)\cap E(x,z)$ is $(n,r)$.   So by compactness there iis $psi_{p}(y)\in p$, and single $\chi(y,z)$ of dimension $< 2k$ such that  $\psi_{p}(y)$ and $\psi_{p}(z)$ and $\neg\chi(y,z)$ implies $\delta(y,z)$.  Now again by compactness finitely many of the $\psi_{p}$ cover $W$ up to an $acl(A)$-definable set of dimension $< k$. 
The upshot is that we can partition $W$ into finitely many $acl(A)$-definable sets $W_{1},..,W_{t}$ say, such that for each $i = 1,..,t$, either for all $(a,b)\in W_{i}\times W_{i}$ except for an $acl(A)$-definable subset of dimension $< 2k$ we have $dim(E(x,a)\wedge E(x,b)) < n$, or for all $(a,b)\in W_{i}\times W_{i}$ except for an $acl(A)$-definable subset of $dim < 2k$ we have $(dim, meas)(E(x,a)\wedge E(x,b)) = (n, r_{i})$ for a fixed positive rational $r_{i}$.

We now have to deal with what goes on for $(a,b)\in  W_{i}\times W_{j}$ for $i\neq j$. For simplicity of presentation we consider the case where $t = 2$. Fix complete $q(z)$ over $acl(A)$ implying $z\in W_{2}$. As above we find a finite partition $P_{q}$ of $W_{1}$ into $acl(A)$-definable sets, and an $acl(A)$-definable subset $W_{2,q}$ of $W_{2}$ which contains $q$, and such that for each $Y\in P_{q}$ the pair $Y$, $W_{2,q}$ is {\em good} in the obvious sense that $dim(E(x,a)\cap E(x,b)) < n$ for almost all $(a,b)\in Y\times W_{2,q}$ or $(dim,meas)(E(x,a)\cap E(x,b)) = (n,r)$ for almost all $(a,.b)\in Y\times W_{2,q}$ (fixed $r$).  By compactness finitely many $W_{2,q}$ cover $W_{2}$ up to $dim < n$. This gives a partition of $W_{2}$. Take the intersection of the partitions $P_{q}$ of $W_{1}$  (for the relevant finite set of $q$'s), and we obtain the required partition of $W$. 

\vspace{5mm}
\noindent
{\em Additional remark.}
\newline
After we wrote this note, Hrushovski sent us his own commentary on  \cite{Tao}, including an essentially identical account of Tao's Proposition 27 to that above. Then Tao himself  gave another proof valid in all characteristics \cite{Tao2}. Nevertheless we think it is worth having the proof above available.

\end{document}